\documentclass[reqno,a4paper,11pt]{amsart}
\usepackage{soul}

\usepackage{amsmath,amssymb,amsthm}

\usepackage{mathrsfs}

\usepackage{mathtools}
\usepackage{commath}

\usepackage{thmtools}
\usepackage{thm-restate}

\usepackage{cases}
\usepackage{enumitem}
\setlist[enumerate,1]{label=(\roman*)}

\usepackage[pdftex, pdfborderstyle={/S/U/W 0}]{hyperref}
\hypersetup{
    colorlinks=true,
    linkcolor=magenta,
    citecolor=cyan,
}

\usepackage{cleveref}

\usepackage{etoolbox}

\usepackage{comment}

\usepackage[numbers]{natbib}

\numberwithin{equation}{section}

\ifdefined\thmcolor
\declaretheoremstyle[
  shaded={bgcolor=\thmcolor}
]{plain}
\else
\fi

\ifdefined\defcolor
\declaretheoremstyle[
  headfont=\normalfont\bfseries,
  bodyfont=\normalfont,
  shaded={bgcolor=\defcolor}
]{noital}
\else
\declaretheoremstyle[
  headfont=\normalfont\bfseries,
  bodyfont=\normalfont,
]{noital}
\fi

\declaretheorem[style=plain,numberwithin=section,name=Theorem]{theorem}
\declaretheorem[style=plain,sibling=theorem,name=Proposition]{proposition}
\declaretheorem[style=plain,sibling=theorem,name=Lemma]{lemma}

\declaretheorem[style=plain,sibling=theorem,name=Conjecture]{conjecture}

\declaretheorem[style=plain,sibling=theorem,name=Question]{question}
\declaretheorem[style=plain,sibling=theorem,name=Observation]{observation}

\declaretheorem[style=plain,numbered=no,name=Theorem]{theorem-n}
\declaretheorem[style=plain,numbered=no,name=Proposition]{proposition-n}
\declaretheorem[style=plain,numbered=no,name=Lemma]{lemma-n}
\declaretheorem[style=plain,numbered=no,name=Corollary]{corollary-n}
\declaretheorem[style=plain,numbered=no,name=Conjecture]{conjecture-n}
\declaretheorem[style=plain,numbered=no,name=Claim]{claim-n}
\declaretheorem[style=plain,numbered=no,name=Fact]{fact-n}
\declaretheorem[style=plain,numbered=no,name=Open Problem]{openproblem-n}
\declaretheorem[style=plain,numbered=no,name=Question]{question-n}
\declaretheorem[style=plain,numbered=no,name=Observation]{observation-n}

\declaretheorem[style=noital,sibling=theorem,name=Definition]{definition}

\declaretheorem[style=noital,sibling=theorem,name=Example]{example}

\declaretheorem[style=noital,numbered=no,name=Remark]{remark-n}
\declaretheorem[style=noital,numbered=no,name=Definition]{definition-n}
\declaretheorem[style=noital,numbered=no,name=Construction]{construction-n}
\declaretheorem[style=noital,numbered=no,name=Example]{example-n}

\newcommand{\defined}{\mathrel{\coloneqq}}

\let\st\relax
\newcommand{\st}{\mathbin{\colon}}

\undef{\set}
\DeclarePairedDelimiter{\set}{\lbrace}{\rbrace}

\undef{\emptyset}
\newcommand{\emptyset}{\varnothing}

\newcommand{\union}{\mathbin{\cup}}
\newcommand{\inter}{\mathbin{\cap}}

\newcommand{\from}{\colon}

\undef{\abs}
\DeclarePairedDelimiterX{\abs}[1]
  {\lvert}{\rvert}{\ifblank{#1}{\,\cdot\,}{#1}}

\undef{\norm}
\DeclarePairedDelimiterX{\norm}[1]
  {\lVert}{\rVert}{\ifblank{#1}{\,\cdot\,}{#1}}

\DeclarePairedDelimiterX{\inner}[2]
  {\langle}{\rangle}{\ifblank{#1}{\,\cdot\,}{#1},\ifblank{#2}{\,\cdot\,}{#2}}

\DeclareMathDelimiter{\given}
  {\mathbin}{symbols}{"6A}{largesymbols}{"0C}

\DeclareMathOperator{\Prob}{\mathbb{P}}
\DeclarePairedDelimiterXPP{\prob}[1]
  {\Prob}{\lparen}{\rparen}{}
  {\renewcommand{\given}{\nonscript\;\delimsize\vert\nonscript\;\mathopen{}}#1}

\DeclareMathOperator{\Expec}{\mathbb{E}}
\DeclarePairedDelimiterXPP{\expec}[1]
  {\Expec}{\lparen}{\rparen}{}
  {\renewcommand{\given}{\nonscript\;\delimsize\vert\nonscript\;\mathopen{}}#1}

\DeclareMathOperator{\Var}{Var}
\DeclarePairedDelimiterXPP{\var}[1]
  {\Var}{\lparen}{\rparen}{}
  {\renewcommand{\given}{\nonscript\;\delimsize\vert\nonscript\;\mathopen{}}#1}

\DeclareMathOperator{\Cov}{Cov}
\DeclarePairedDelimiterXPP{\cov}[2]
  {\Cov}{\lparen}{\rparen}{}{#1,#2}

\newcommand{\sseq}{\subseteq}

\let\l\relax
\newcommand{\l}{\ell}

\newcommand{\NN}{\mathbb{N}}

\newcommand{\QQ}{\mathbb{Q}}

\newcommand{\cF}{\mathcal{F}}

\newcommand{\gk}{\kappa}
\newcommand{\gl}{\lambda}

\newcommand{\gw}{\omega}

\usepackage{geometry}
\usepackage{centernot}

\usepackage{tikz}
\usetikzlibrary{quotes} 
\usetikzlibrary{calc}

\let\th\relax
\newcommand{\th}{\textsuperscript{th}}

\makeatletter
\newcommand{\oset}[3][0ex]{%
  \mathrel{\mathop{#3}\limits^{
    \vbox to#1{\kern-2\ex@
    \hbox{$\scriptstyle#2$}\vss}}}}
\makeatother

\newcommand{\incomp}{\parallel}
\newcommand{\comp}{\sim}

\let\and\relax
\newcommand{\and}{\wedge}

\begin{document}

\title{The Aharoni--Korman conjecture is false}

\author{Lawrence Hollom}
\address{Department of Pure Mathematics and Mathematical Statistics (DPMMS), University of Cambridge, Wilberforce Road, Cambridge, CB3 0WA, United Kingdom}
\email{lh569@cam.ac.uk}



\begin{abstract}
  A poset $P$ is said to satisfy the \emph{finite antichain condition}, or \emph{FAC}, if it has no infinite antichain.
  It was conjectured by Aharoni and Korman in 1992 that any FAC poset $P$ possesses a chain $C$ and a partition into antichains such that $C$ meets every antichain of the partition.
  In this work we provide a counterexample to this conjecture, demonstrating that it is false.
  We also discuss variations of the conjecture which may yet be true.
\end{abstract}

\maketitle


\section{Introduction}
\label{sec:intro}

Posets with no infinite antichains have a rich and intricate structure, but are relatively under-studied as objects of interest in their own right.
One of the most significant problems in this area is due to Aharoni and Korman, who conjectured in 1992 \cite{aharoni1992greene} that any poset $P$ with no infinite antichain possesses a chain $C$ and a partition into antichains such that $C$ meets every antichain of the partition.
The problem arose naturally from -- and has strong links to -- the study of covers and matchings in infinite hypergraphs.

This conjecture is often referred to as ``the fishbone conjecture'' (due to how the structure of the chains and antichains resembles the skeleton of a fish), and there have been many attempts to prove the conjecture; indeed, it has been described as ``one of the most attractive'' conjectures in the area \cite{aharoni2009menger}.
However, despite this attention, little progress has been made on this problem.

We demonstrate why the conjecture has proven so resilient, and provide a counterexample to the conjecture.

The conjecture of Aharoni and Korman can be stated formally and in full generality as follows.

\begin{conjecture}[\cite{aharoni1992greene}, Conjecture 4.1]
    \label{conj:ak}
    If a poset $P$ contains no infinite antichain then, for every positive integer $k$, there exist $k$ chains $C_1,\dotsc,C_k$ and a partition of $P$ into disjoint antichains $(A_i\st i\in I)$ such that each $A_i$ meets $\min(\abs{A_i},k)$ chains $C_j$.
\end{conjecture}

This conjecture -- or sometimes the $k=1$ case thereof -- is referred to as \emph{the Aharoni--Korman conjecture} or \emph{the fishbone conjecture}.
We will use the former name, and work exclusively with the case of $k=1$.

To give an example, if $P$ is the poset on the set $\NN\times\NN$ with ordering given by setting $(x,y)\leq(u,v)$ if and only if $x\leq u$ and $y\leq v$, then the $k=1$ case of \Cref{conj:ak} holds by taking $C = \set{(0,y)\st y\in\NN}$ and $A_i=\set{(x,y)\in P\st x+y=i}$ for all integers $i\geq 0$.
However, this partition into antichains is not unique. For example, we could instead take all antichains of the form $\set{(x,y),((x+y)(x+y+1)/2 + y, 0)}$ for $x,y \geq 0$, and these are also pairwise disjoint, partition $P$, and all meet $C$.

Our main result is the following.

\begin{theorem}
    \label{thm:counterexample}
    There is a countable poset $P$ for which \Cref{conj:ak} is false.
\end{theorem}

As is discussed further in \cite{hollom2024resolution}, the poset $P$ constructed in the proof of \Cref{thm:counterexample} is necessarily quite complex, as we now explain.

We say that a chain $C$ in a poset $P$ is \emph{contiguous} if there is no $x\in P\setminus C$ and $y,z\in C$ with $y<x<z$ and $C\union\set{x}$ a chain.
The \emph{linear sum} $\bigoplus_{i\in I}P_i$ of a collection $(P_i\st i\in I)$ of posets indexed by a total order $I$ is the ordering on $\bigsqcup_{i\in I}P_i$ wherein all elements of $P_j$ are set to be greater than all elements of $P_i$ whenever $j > i$ in $I$.
Following standard set-theoretic notation, we write $\gw$ or $\NN$ for the set of natural numbers. In particular, we use $\NN$ to refer to the (unordered) set of natural numbers, and $\gw$ to refer to the poset naturals with the standard ordering.
We write $C^*$ for the poset equal to $C$ with the ordering reversed.
With this notation in hand, we can state the following result of \cite{hollom2024resolution}.

\begin{theorem}[\cite{hollom2024resolution}]
    \label{thm:main-on-first-page}
    Let $P$ be a countable poset with no infinite antichain such that, for any contiguous chain $C\sseq P$, neither $C$ nor $C^*$ can be written as $\bigoplus_{i\in\gw} C_i$ with each $C_i$ infinite and co-wellfounded.
    Then there is a chain $C\sseq P$ and partition $P=\bigcup_{i\in I} A_i$ into antichains such that $C$ meets every antichain $A_i$.
\end{theorem}

We thus see that a counterexample to \Cref{conj:ak} must contain a wellfounded linear sum of co-wellfounded posets, and so necessarily has a somewhat complex structure. 
Indeed, the poset we will construct in \Cref{ex:counterexample} can be seen to contain such a sum.


\subsection{Terminology and notation}
\label{subsec:basics}

We now recall some fundamental facts and definitions that we will make use of throughout the rest of this paper.

A subset $X$ of a partially ordered set, or poset, $P$, is a \emph{chain} if the elements of $X$ are pairwise comparable, and it is an \emph{antichain} if its elements are pairwise incomparable.
If $P$ has no infinite antichain, then we say that it satisfies the \emph{finite antichain condition}, and for brevity we will refer to $P$ as an \emph{FAC poset}.
We note here that the \emph{width} of an FAC poset $P$ -- the supremum of the sizes of antichains in $P$ -- may in fact be infinite.

Given two posets $P$ and $Q$ with orders $\leq_P$ and $\leq_Q$ respectively, the \emph{Cartesian product} poset $P\times Q$ is is the poset on the Cartesian product set $P\times Q$, with order $\leq$ given by $(p,q)\leq (p',q')$ if and only if $p\leq_P p'$ and $q\leq_Q q'$.

For an arbitrary poset $P$, a subset $X\sseq P$ is said to be \emph{convex}, or an \emph{interval} of $P$ if for all $x,y\in X$, if $z\in P$ has $x<z<y$, then $z\in X$.

We use the symbol $\comp$ to represent comparability in a poset.
That is, if $x,y\in P$ have $x\leq y$ or $y\leq x$, then we write $x\comp y$.
Similarly, we write $X\comp Y$ for $X,Y\sseq P$ if $x\comp y$ for every $x\in X$ and $y\in Y$.

Finally, we provide the following definition for ease of discussion.

\begin{definition}
    \label{def:spine}
    Given a poset $P$, we say that a subset $C\sseq P$ is a \emph{spine} of $P$ if $C$ is a chain and there is a partition of $P$ into antichains $(A_i\st i\in I)$ such that $C$ meets every antichain $A_i$.
\end{definition}


\subsection{Known results}
\label{subsec:lit-review}

Questions concerning the partition of posets into various structures can be traced back to Dilworth's theorem \cite{dilworth1950decomposition}, originally published in 1950, which states that a poset $P$ with width at most $k$ has a partition into at most $k$ chains.
As will become a theme throughout this introduction, the generalisation to infinite posets presented some further challenges.
Indeed, Perles \cite{perles1963dilworth} noted that the strongest ``infinite Dilworth's theorem'' one could hope for -- that any poset of width at most $\gl$ can be partitioned into at most $\gl$-many chains -- is false, as, for example, for any infinite cardinal $\gk$, the poset $\gk\times\gk$ has no infinite antichain, and thus width $\aleph_0$, but cannot be partitioned into fewer than $\gk$-many chains.
However, Dilworth's theorem can be generalised to infinite FAC posets $P$ under the assumption that $P$ has some partition into finitely many chains, and this was proved by Abraham in 1987 \cite{abraham1987note}.

In 1976, Greene and Kleitman \cite{greene1976structure} generalised Dilworth's theorem to families of antichains, proving that, for any finite poset $P$ and integer $k$, the maximal size of the union of $k$ antichains of $P$ is equal to the minimal value of $\sum_{i\leq m}\min(\abs{C_i},k)$ over all partitions $\set{C_1,\dotsc,C_m}$ of $P$ into chains.
The Greene--Kleitman theorem was then generalised to the infinite setting by Aharoni and Korman \cite{aharoni1992greene}.
More precisely, they proved the following theorem, which they described as ``the `correct' infinite version of Greene-Kleitman’s theorem.''

\begin{theorem}[\cite{aharoni1992greene}, Theorem 3.1]
    \label{thm:infinite-greene-kleitman}
    Let $P$ be a poset with no infinite chains and let $k$ be a positive integer. There exists then a partition $(C_i\st i\in I)$ of $P$ into disjoint chains and $k$ disjoint antichains $A_j$ $(1 \leq j \leq k)$ such that every chain $C_i$ meets $\min(k, \abs{C_i})$ antichains $A_j$.
\end{theorem}

Having proved \Cref{thm:infinite-greene-kleitman}, Aharoni and Korman then asked about the dual version, interchanging the roles of chains and antichains, which lead them to \Cref{conj:ak}.
This conjecture was later referred to by Aharoni as the ``fishbone conjecture'' (see for example \cite{aharoni2022strongly}), but is often called ``the Aharoni--Korman conjecture''.
The inspiration for the conjecture in fact came from the consideration of strongly minimal covers and strongly maximal matchings in hypergraphs (see for example \cite{aharoni1991infinite}), and has close links to infinite graph theory as a whole.
One such related area is the generalisation of Menger's theorem to infinite graphs.
This problems has a rich history, and led to many partial results \cite{,aharoni1987menger,aharoni1994menger,diestel2003countable} before the seminal paper of Aharoni and Berger \cite{aharoni2009menger} proved the full generalisation.
Indeed, in \cite{aharoni2009menger}, \Cref{conj:ak} was described as ``one of the most attractive'' conjectures in the area of infinite matching theory.

Aharoni and Korman also posed several other, more general conjectures in \cite{aharoni1992greene}.
In 1995, Aharoni and Loebl \cite{aharoni1995strongly} resolved one of these conjectures, concerning perfect graphs, in some special cases.
Much more recently, the strongest of the conjectures from \cite{aharoni1992greene}, concerning strongly minimal covers of hypergraphs, has been disproved by van der Zypen \cite{van2022counterexample}.

While there have been some positive results concerning \Cref{conj:ak}, the previous state of the art was still a considerable way from a resolution of the conjecture.
Indeed, the case of $k=1$ and posets of width at most $2$ was dealt with by Aharoni and Korman themselves \cite{aharoni1992greene} by application of an infinite version of K\"{o}nig's theorem due to Aharoni \cite{aharoni1984konig}.
The first partial result after the conjecture was introduced was due to Duffus and Goddard \cite{duffus2002intervals} who proved that \Cref{conj:ak} holds in two different cases.
Firstly, they proved the conjecture for posets of the form $C\times P$, where $C$ is a chain and $P$ is finite (this result in fact first appeared in Goddard's doctoral thesis \cite{goddard1996ordered}).
Secondly, they proved the conjecture for FAC posets with no infinite intervals.

Since then, Zaguia has recently proved a series of results, now collected into a single paper \cite{zaguia2024progress}, showing that \Cref{conj:ak} is true when the FAC poset $P$ is $N$-free (that is, if $w,x,y,z\in P$ have $w\leq y$, $x\leq y$, and $x\leq z$, then some further comparison must hold between these four points), or has locally finite incomparability graph.
The \emph{incomparability graph} of a poset $P$ is a graph $G$ with vertex set $P$, and $xy$ is an edge of $G$ if and only if $x\incomp y$ in $P$.

Recently, the author proved that any FAC poset satisfying a particular condition on decompositions of its chains must satisfy the Aharoni--Korman conjecture.
The following definition is central to this result.

\begin{definition}
    \label{def:vacillating}
    We say that a poset $P$ is \emph{vacillating} if there is no contiguous chain $C$ in $P$ for which either $C$ or its reverse $C^*$ can be written as $\bigoplus_{i\in \gw} C_i$, where each $C_i$ is infinite and co-wellfounded.
\end{definition}

With this definition in hand, one of the most significant positive partial results concerning the Aharoni--Korman conjecture is \Cref{thm:main-on-first-page}, which can now be succinctly stated as the fact that any vacillating FAC poset has a spine.

Taking a more broad view, beyond just partial results on \Cref{conj:ak}, FAC posets in general have a rich structure, but have received only sparse attention.
Now that we know that FAC posets do not necessarily contain spines, one may ask how \Cref{conj:ak} could be weakened so that it may become true.
One possible direction would be in the study of strongly maximal structures.

Given a set family $\cF$, a set $A\in \cF$ is \emph{maximal} if $B\supseteq A$ and $B\in \cF$ implies that $A=B$.
Moreover, $A$ is \emph{strongly maximal} if, for every $B\in \cF$, we have $\abs{A\setminus B} \geq \abs{B\setminus A}$.
As was noted and discussed by Aharoni \cite{aharoni1991infinite}, it is a consequence of \Cref{thm:infinite-greene-kleitman} that if $P$ is a poset with no infinite chain, then $P$ contains a strongly maximal antichain.
Later, Aharoni and Berger \cite{aharoni2011strongly} strengthened this result by producing a weaker condition which also implies the existence of a strongly maximal antichain.

In a direction more analogous with spines in posets, one may consider strongly maximal chains as a weaker form of spines, which we record the definition of as follows.

\begin{definition}
    \label{def:smc}
    A chain $C$ of a poset $P$ is a \emph{maximal chain} if no element of $P\setminus C$ is comparable to every element of $C$.
    $C$ is a \emph{strongly maximal chain} if for every chain $D\sseq P$, we have that $\abs{C\setminus D}\geq \abs{D\setminus C}$.    
\end{definition}

Indeed, strongly maximal chains were the subject of an investigation of the author in \cite{hollom2024resolution}, and it is therein proved that countable FAC posets do indeed contain strongly maximal chains.
We record the following observation, the proof of which is simple, and not given here.
\begin{observation}
    \label{obs:spine-is-smc}
    If $C$ is a spine of a poset $P$, then $C$ is also a strongly maximal chain in $P$.
\end{observation}
The above observation allows us to consider the existence of strongly maximal chains as a resolution of a weaker version of the Aharoni--Korman conjecture.

\section{Proof of main result}

We now present a counterexample to the Aharoni--Korman conjecture, \Cref{conj:ak}.
Due to the results of \cite{hollom2024resolution} discussed in \Cref{sec:intro}, we know that if a poset $P$ is to be a counterexample to \Cref{conj:ak}, then $P$ cannot be vacillating. 
Indeed, our counterexample will consist of an infinite non-increasing sequence of copies of $\gw\times\gw$, along with some extra structure.
We now define this counterexample.

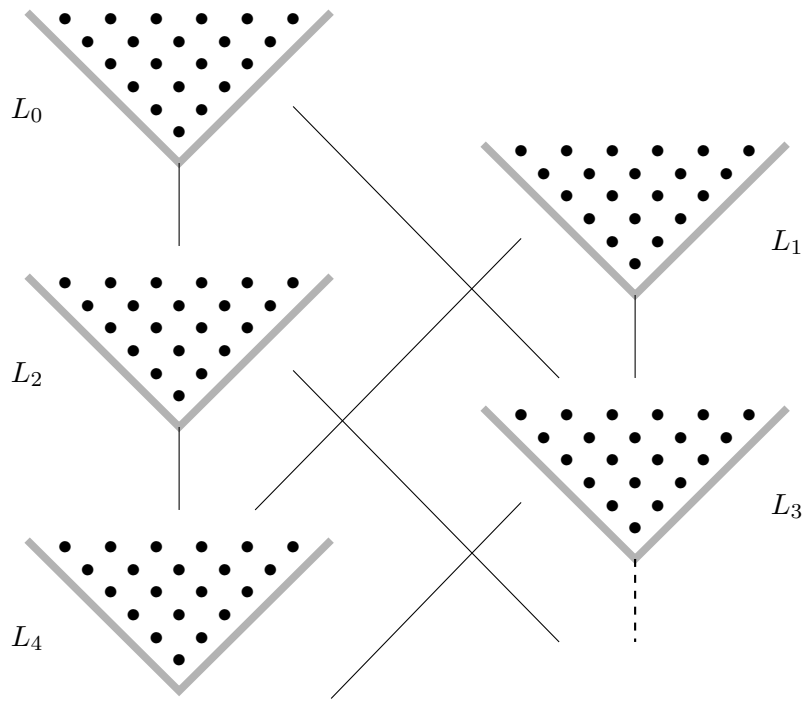
\begin{figure}[htbp]
    \centering
    \begin{tikzpicture}[scale=1]
    \pgfmathsetmacro\spacing{0.3}
    \pgfmathsetmacro\trianglesize{2}
    \pgfmathsetmacro\yoff{0.2}

    \definecolor{grey}{rgb}{0.7,0.7,0.7}

    \foreach \xbase/\ybase [count = \i from 1] in {0/0, 6/1.75, 0/3.5, 6/5.25, 0/7} {
        \coordinate (l1\i) at (\xbase, \ybase - \yoff);
        \coordinate (l2\i) at (\xbase - \trianglesize, \ybase + \trianglesize - \yoff);
        \coordinate (l3\i) at (\xbase + \trianglesize, \ybase + \trianglesize - \yoff);
        \coordinate (top\i) at (\xbase, \ybase + \trianglesize + \yoff);
        \draw[color=grey,line width = 1mm] (l2\i) -- (l1\i) -- (l3\i);

        \foreach \s in {0,...,5} {
            \foreach \x in {0,...,\s} {
                \pgfmathsetmacro\xoff{(2 * \x - \s) * \spacing}
                \node (v\i\s\x) at (\xbase + \xoff,\ybase + \s * \spacing + \yoff) {$\bullet$};
            }
        }
    }

    \coordinate (top0) at (6, 0.45);

    \draw (top1) -- (l13);
    \draw (top2) -- (l14);
    \draw (top3) -- (l15);
    \draw[thick,dashed] (l12) -- (top0);

    \draw (1,2.2) -- (4.5, 5.8);
    \draw (2,-0.3) -- (4.5, 2.3);
    \draw (5,3.95) -- (1.5, 7.55);
    \draw (5,0.45) -- (1.5, 4.05);

    \node (L0) at (-2,7.5) {$L_0$};
    \node (L2) at (-2,4) {$L_2$};
    \node (L4) at (-2,0.5) {$L_4$};

    \node (L1) at (8,5.75) {$L_1$};
    \node (L3) at (8,2.25) {$L_3$};

\end{tikzpicture}
    \caption{An approximate Hasse diagram of $P$. The regions within the grey V-shapes are order-isomorphic to $\gw\times\gw$, and this structure continues infinitely downwards. Relations between $L_n$ and $L_{n+1}$ are shown in more detail in \Cref{fig:counterexample-details}.}
    \label{fig:counterexample}
\end{figure}

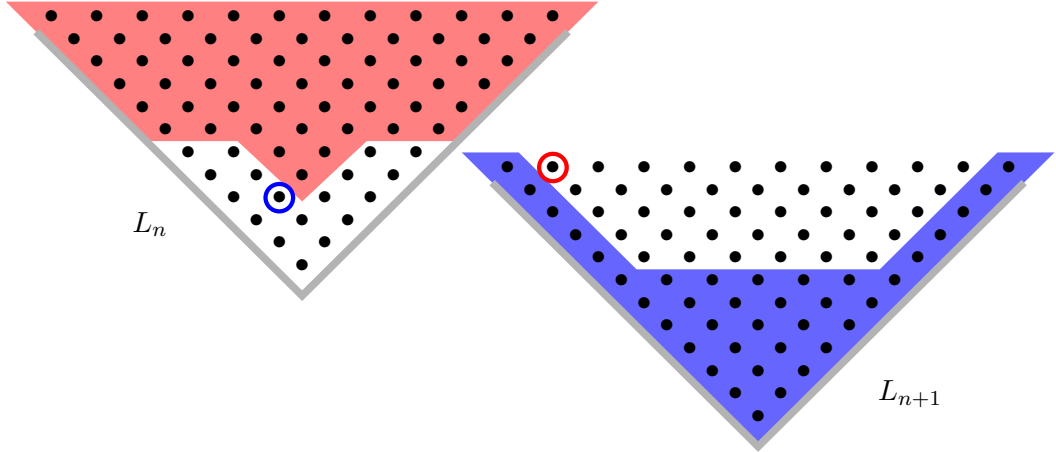
\begin{figure}[ht]
    \centering
    \begin{tikzpicture}[scale=1]
    \pgfmathsetmacro\spacing{0.3}
    \pgfmathsetmacro\trianglesize{3.5}
    \pgfmathsetmacro\yoff{0.2}

    \definecolor{grey}{rgb}{0.7,0.7,0.7}

    \fill[blue, opacity=0.6] (6,-0.2) -- (2.1,3.7) -- (2.85,3.7) -- (4.4,2.15) -- (7.6,2.15) -- (9.15,3.7) -- (9.9,3.7) -- cycle;
    \fill[red, opacity=0.5] (-2.05,3.85) -- (-3.9,5.7) -- (3.9,5.7) -- (2.05,3.85) -- (0.85,3.85) -- (0,3.05) -- (-0.85,3.85) -- cycle;

    \foreach \xbase/\ybase [count = \i from 1] in {6/0, 0/2} {
        \coordinate (l1\i) at (\xbase, \ybase - \yoff);
        \coordinate (l2\i) at (\xbase - \trianglesize, \ybase + \trianglesize - \yoff);
        \coordinate (l3\i) at (\xbase + \trianglesize, \ybase + \trianglesize - \yoff);
        \coordinate (top\i) at (\xbase, \ybase + \trianglesize + \yoff);
        \draw[color=grey,line width = 1mm] (l2\i) -- (l1\i) -- (l3\i);

        \foreach \s in {0,...,11} {
            \foreach \x in {0,...,\s} {
                \pgfmathsetmacro\xoff{(2 * \x - \s) * \spacing}
                \node (v\i\s\x) at (\xbase + \xoff,\ybase + \s * \spacing + \yoff) {$\bullet$};
            }
        }
    }

    \node[circle,color=blue,draw,ultra thick] (upbright) at (- \spacing, 2 + 3 * \spacing + \yoff) {};
    \node[circle,color=red,draw,ultra thick] (downbright) at (6 - 9 * \spacing, 11 * \spacing + \yoff) {};

    \node (L0) at (-2,2.75) {$L_n$};
    \node (L1) at (8,0.5) {$L_{n+1}$};

\end{tikzpicture}
    \caption{Two examples of the relations between $L_n$ and $L_{n+1}$. The red region in $L_n$ is the set of those points above $(1,10,n+1)$ (shown in a red circle), and the blue region in $L_{n+1}$ is the set of those points below $(1,2,n)$ (shown in a blue circle).}
    \label{fig:counterexample-details}
\end{figure}

\begin{example}
    \label{ex:counterexample}
    Let $P$ be the poset defined on the set $\NN \times \NN \times \NN$ with the order $\leq$ on elements $(x,y,n)$ of $P$ given by setting $(x,y,n)\leq (u,v,m)$ if and only if
    \begin{itemize}
        \item $n\geq m+2$, or
        \item $n = m$ and $x \leq u$ and $y \leq v$, or
        \item $n = m+1$ and $\min\set{x,y} + 1 \leq \min\set{u,v}$, or
        \item $n = m+1$ and $x+y \leq 2(u+v)$.
    \end{itemize}
    We note without proof that the above relation is in fact a partial order.
    Define the set $L_n\sseq P$ to be the \emph{$n$\th\ level of $P$}, defined as $L_n\defined \set{(x,y,n)\st x,y\in\gw}$.
    A Hasse diagram of the high-level structure of $P$ is shown in \Cref{fig:counterexample}, and two examples showing the relations between $L_n$ and $L_{n+1}$ are shown in \Cref{fig:counterexample-details}.
\end{example}

We remark before continuing that the final condition in the above definition would work just as well with the number $2$ replaced by any real number $\alpha > 1$, we merely use 2 for convenience. 
Moreover, this condition is only used at the very end of the proof.

We now prove the following proposition, which states that $P$ is a counterexample to \Cref{conj:ak}.

\begin{proposition}
    \label{prop:spineless}
    $P$ is a countable FAC poset, but there is no chain $C\sseq P$ which is a spine of $P$. 
\end{proposition}

Before proceeding to the proof, we remark that our proof of \Cref{prop:spineless} has been formally verified using the Lean theorem prover by Mehta, and the code is available online as part of the mathlib4 library \cite{mehta2025formal}.

We know that, due to \Cref{thm:main-on-first-page}, a counterexample to \Cref{conj:ak} cannot be vacillating.
Indeed, the set $\set{(x,0,2n)\st x,n\in \NN} \sseq P$ is order-isomorphic to $\bigoplus_{i\in \gw^*} \gw$, and thus $P$ is not vacillating.

We will prove \Cref{prop:spineless} via a series of observations and lemmas.
Fix for the rest of this section a chain $C\sseq P$, which we will prove is not a spine of $P$.
For ease of reference, we first enumerate a series of observations.

\begin{observation}
    \label{obs:obs}
    The following facts are all immediate.
    \begin{enumerate}
        \item For all $n$, $L_n$ is an interval of $P$.
        \item \label{obs:antichains} For any $n$ and $s$, the set $K_{n,s}\defined\set{(x,y,n)\st x+y=s}\sseq L_n$ is an antichain of $P$.
        \item \label{obs:chains} For any fixed $n$, the induced order on $L_n\sseq P$ is isomorphic to $\gw\times\gw$.
        \item \label{obs:interpolate} For any three points $(x,y,n)\leq (u,v,n) \leq (w,z,n)$ in $L_n$, there is a chain $D\sseq [(x,y,n),(w,z,n)]$ containing $(u,v,n)$ of length $w+z+1-x-y$.
    \end{enumerate}
\end{observation}

\begin{lemma}
    \label{lem:fac}
    $P$ contains no infinite antichain.
\end{lemma}

\begin{proof}  
    Assume for contradiction that there is some infinite antichain $A\sseq P$.
    As $L_n\comp L_m$ whenever $\abs{n-m}\geq 2$, we know that $A$ can intersect at most two levels of $P$, and thus there must be some level $L_n$ such that $A\inter L_n$ is infinite.
    We may thus assume that in fact $A\sseq L_n$ for some integer $n$.
    
    However, by \Cref{obs:obs} \cref{obs:chains}, $L_n$ is order-isomorphic to $\omega \times \omega$, and it is easy to check that $\gw\times\gw$ is an FAC poset.
    Thus $P$ has no infinite antichains, as required.
\end{proof}

Define a sequence $m_0,m_1,\dotsc$ of integers by
\begin{align*}
    m_n\defined \min\set{\min\set{x,y}\st (x,y,n)\in C}.
\end{align*}

\begin{lemma}
    \label{lem:exists-sparse-level}
    There is some integer $n\geq 1$ such that the size of the intersection $\abs{L_n\inter C}$ of $C$ with the $n$\th\ level $L_n$ of $P$ is finite.
\end{lemma}

\begin{proof}
    Assume for contradiction that $\abs{L_n\inter C}$ is infinite for every integer $n$.
    We claim that $m_1>m_2>\cdots$, which will immediately yield a contradiction.
    Indeed, fix some $n$, and let $(u,v,n)\in C$ have $\min\set{u,v}=m_n$.
    Consider $L_{n+1}\inter C$; this set is infinite, and so due to \Cref{obs:obs} \cref{obs:antichains}, there must be a infinitely many points $(x,y,n+1)\in C$ such that $x+y > 2(u+v)$.
    Thus to have $(x,y,n+1) < (u,v,n)$, we must have $m_{n+1}\leq \min\set{x,y}+1\leq \min\set{u,v} = m_n$.
    Therefore $m_{n+1} < m_n$ for all $n$, leading to a contradiction, as required.
\end{proof}

Fix $n$ as in \Cref{lem:exists-sparse-level} so that $\abs{L_n\inter C}$ is finite; we will prove that there is no function $f\from C\union L_n \union L_{n+1}\to C$ such that $f$ acts as the identity on $C$, and the preimage of any point $x\in C$ is an antichain.
Note that the existence of such a function is equivalent to the existence of a partition into antichains (the sets in the partition are simply the preimages of $f$), and that if there does not exist such a function for $C\union L_n\union L_{n+1}$, then there certainly cannot exist such a function for the whole of $P$.

Assume for contradiction that some such function does exist, and thus fix for the rest of this section a function $f\from C\union L_n \union L_{n+1}\to C$ which acts as the identity on $C$ and for which every preimage is an antichain of $P$.

We first show that there are infinitely many points in $L_n$ which $f$ does not send to $C\inter L_n$ or $C\inter L_{n+1}$.
Indeed, define $R\sseq L_n$ to be those elements which are greater than all elements of $C\inter L_n$:
\begin{align*}
    R\defined\set{q\in L_n\st q > \max(C\inter L_n)}.
\end{align*}
As $C\inter L_n$ is finite due to \Cref{lem:exists-sparse-level}, we know that there is some integer $a$ such that $\set{(x,y,n)\st x\geq a \and y\geq a}\sseq R$.
As $R$ is entirely comparable with $C\inter L_n$ by definition, $f$ cannot send any element of $R$ to $C\inter L_n$, and so must send these elements to $C\inter L_{n-1}$ or $C\inter L_{n+1}$.

If $C\inter L_{n+1}$ is non-empty then, noting that $L_{n+1}$ is well-founded, and so any chain it contains is well-ordered, set $(x_0,y_0,n+1) = \min(C\inter L_{n+1})$, note that $m_{n+1} = \min\set{x_0,y_0}$, and define a set $S\sseq R$ as follows.
\begin{align*}
    S\defined
    \begin{cases}
        \set{(x,y,n)\in R\st (x,y,n)\comp C\inter L_{n+1}} &\text{if } C\inter L_{n+1}\text{ is finite,}\\
        \set{(x,y,n)\in R\st \min\set{x,y}\geq \max\set{x_0,y_0} + 1} &\text{otherwise.}
    \end{cases}
\end{align*}

Note that $S$ is always infinite, and that if $(x_0,y_0,n+1)$ does not exist, then we are in the first of the above cases, so do not attempt to make use of its value.
Moreover, as was also the case for $R$, note that there is some integer $b$ such that $\set{(x,y,n)\st x\geq b \and y\geq b}\sseq S$.
We now prove that $f$ cannot send any element of $S$ to $C\inter L_{n+1}$.

\begin{lemma}
    \label{lem:dont-map-down}
    If $p=(x,y,n)\in S$, then $f(p)\notin C\inter L_{n+1}$.
\end{lemma}

\begin{proof}
    If $C\inter L_{n+1}$ is finite, then the result is immediate as $S$ is entirely comparable to $C\inter L_{n+1}$ by definition.
    We may thus assume that $C\inter L_{n+1}$ is infinite.
    Define the chains
    \begin{align*}
        D_0\defined \set{(x,y_0,n+1)\st x\geq x_0} \quad\text{and}\quad D_1\defined \set{(x_0,y,n+1)\st y\geq y_0}.
    \end{align*}
    Note that at least one of the chains $D_0$ and $D_1$ has the property that each of its elements is below some element of $C\inter L_{n+1}$. Assume without loss of generality (by exchanging the roles of $x$ and $y$ if necessary) that this holds for $D_0$.

    Recall the definition of $K_{n,s}$ from \Cref{obs:obs} \cref{obs:antichains}.
    We now claim that $\abs{C\inter K_{n+1,s}} = 1$ for every $s\geq x_0+y_0$.
    Firstly, we know from \Cref{obs:obs} that $K_{n+1,s}$ is an antichain, so certainly $\abs{C\inter K_{n+1,s}}\leq 1$.
    If there is some $s_1>x_0 + y_0$ such that $C\inter K_{n+1,s_1} = \emptyset$, then let 
    \begin{align*}
        (x_1,y_1,n+1)\defined\min\set{(x,y,n+1)\in C\st x+y > s_1}.
    \end{align*}
    We thus know that the set $[(x_0,y_0,n+1),(x_1,y_1,n+1)]\inter C$ has size at most $x_1+y_1 - x_0 -y_0$.
    However, as $x_1\geq x_0$ and $y_1\geq y_0$ (as $(x_1,y_1,n+1)\geq (x_0,y_0,n+1)$ in $P$), \Cref{obs:obs} \cref{obs:interpolate} tells us that there is a chain of length $x_1+y_1+1 - x_0 - y_0$ contained in the interval $[(x_0,y_0,n+1),(x_1,y_1,n+1)]$ of $L_{n+1}$.
    Thus we can replace a finite contiguous sub-chain of $C$ with a longer chain, and so $C$ is not a strongly maximal chain.
    However, we know from \Cref{obs:spine-is-smc} that a spine must in particular be a strongly maximal chain, and so this contradicts $C$ being a spine.
    Thus $\abs{C\inter K_{n+1,s}} = 1$ for all $s\geq x_0+y_0$.

    Let $s_0 \defined x_0 + y_0$, and for $s\geq s_0$, let $p_s$ be the unique element of $C\inter K_{n+1,s}$.
    We now claim that 
    \begin{align}
        \label{eq:bijection}
        \forall i\geq 0, \; f((x_0+i,y_0,n+1)) = p_{s_0+i}.
    \end{align}

    In fact, we will prove the following stronger statement, from which \eqref{eq:bijection} follows immediately.
    \begin{align}
        \label{eq:constant-on-rows}
        \forall s \geq s_0 \; \forall \l \geq s-s_0, \; f \text{ is constant on } K_{n+1,s}\inter [p_{s_0}, p_{s_0+\l}].
    \end{align}
    Intuitively, \eqref{eq:constant-on-rows} says that $f$ is ``constant on the rows'' of a suitably defined region of $L_{n+1}$.

    Indeed, if we let $B_\l\defined [p_{s_0},p_{s_0+\l}]$, then it suffices to prove that if $i,j\geq 0$ are such that the points $q\defined (x_0+i+1,y_0+j,n+1)$ and $q'\defined (x_0+i,y_0+j+1,n+1)$ are both in $B_\l$, then we have that $f(q)=f(q')$.

    Note that $q_-\defined (x_0+i,y_0+j,n+1)$ and $q_+\defined (x_0+i+1,y_0+j+1,n+1)$ are both in $B_\l$.
    Thus, applying \Cref{obs:obs} \cref{obs:interpolate}, we see that there are chains $E$ and $E'$ in $B_\l$, both containing $p_{s_0}$, $q_-$, $q_+$, and $p_{s_0+\l}$, with $\abs{E}=\abs{E'}=\l+1$, and differing only in that $E\setminus E' = \set{q}$ and $E'\setminus E = \set{q'}$.

    As $f$ must send all of $B_\l$ to $C\inter B_\l$ (as $B_\l$ is entirely comparable to the rest of $C$), and $\abs{C\inter B_\l} = \l + 1$, we see that $f$ must act as a bijection when restricted to either $E$ or $E'$; we may thus conclude that $f(q)=f(q')$, as required.
    Hence \eqref{eq:constant-on-rows} -- and thus also \eqref{eq:bijection} -- holds.

    Consider a point $p=(x,y,n) \in S$, as in the statement of \Cref{lem:dont-map-down}.
    As $p\in S$, we know that $\min\set{x,y}\geq \max\set{x_0,y_0} + 1$, and so, for any integer $i\geq 0$, $p\geq (x_0+i,y_0,n+1)$.
    But then \eqref{eq:bijection} implies that $p$ is comparable with some element of $f^{-1}(p_{s_0 + i})$ for every $i\geq 0$, and so $f(p)\notin C\inter L_{n+1}$, as required.
\end{proof}

The only remaining possibility is that $f$ sends $S$ to $C\inter L_{n-1}$.
However, $(x,y,n)\in S$ is comparable to all points $(u,v,n-1)\in L_{n-1}$ with $2(u+v)\geq x+y$.
Fix some integer $a\geq 1$ such that the point $(a,a,n)\in S$, and consider the chain
\begin{align*}
    F\defined \set{(a,a,n),(a+1,a,n),\dotsc,(3a-1,a,n),(3a,a,n)}.
\end{align*}
The chain $F$ has length $2a+1$, and is entirely comparable to all points of the form $(u,v,n-1)$ with $u+v\geq 2a$.
Thus $f$ must send the chain $F$ to the intersection of $C$ with the set
\begin{align*}
    T\defined \set{q\in L_{n-1}\st q\not\geq (3a,a,b)} = \set{(u,v,n-1)\st u+v\leq 2a-1}.
\end{align*}
However, the longest chain in $T$ has length $2a$, and so $f$ cannot injectively map $F$ into $C\inter T$, a contradiction.
Therefore the function $f$ cannot exist, and so $P$ has no spine, completing the proof of \Cref{prop:spineless}.
\qed

\section{Concluding remarks and open problems}
\label{sec:conclusion}

In this work we have provided a counterexample to the Aharoni--Korman conjecture.
However, even alongside other work on this conjecture, there are still more questions to ask on this topic.
Indeed, we know from \Cref{thm:main-on-first-page} that any vacillating poset does satisfy the Aharoni--Korman conjecture.
However, there are still families of non-vacillating posets for which the conjecture may hold.
For example, the counterexample presented in \Cref{ex:counterexample} has width $\aleph_0$, and it does not seem easy to adjust the example to have finite width.
For this reason, we are hopeful that the answer to the following question might be positive.

\begin{question}
    \label{q:n-wide}
    Let $P$ be a poset of some finite width $w$.
    Must $P$ have a spine?
\end{question}

We note that it is known (and was proved by Aharoni and Korman themselves \cite{aharoni1992greene}) that every poset of width 2 has a spine.
However, to the best of the author's knowledge, \Cref{q:n-wide} remains open even for countable posets of width 3.

In a slightly different direction, we may weaken the notion of a spine to that of a strongly maximal chain, as defined in \Cref{def:smc}.
While we have shown that FAC posets need not have spines, we conjecture that they do have strongly maximal chains, as follows.

\begin{conjecture}
    \label{conj:big-smcs}
    Let $P$ be an FAC poset.
    Then $P$ has a strongly maximal chain.
\end{conjecture}

It seems that a reasonable starting point may be to try and prove that \Cref{conj:big-smcs} holds for scattered posets (that is, posets which do not admit an embedding of $\QQ$), as this allows for inductive techniques to be applied (see \cite{hausdorff1908classification}).
Finally, it is known that \Cref{conj:big-smcs} holds for countable posets (see \cite{hollom2024resolution}), and so only the case of uncountable posets remains open.

\section{Acknowledgements}

The author would like to thank B\'{e}la Bollob\'{a}s for his thorough reading of the manuscript and many valuable comments. 
The author would also like to thank Bhavik Mehta for producing a formal verification of \Cref{prop:spineless} (available at \cite{mehta2025formal}), and for helping find and resolve several inaccuracies in the original proof of this statement.
Thanks are due to Nikolai Beluhov and George Bergman for each pointing out several inaccuracies in the manuscript, and for further suggestions which significantly improved the presentation of the paper.
Finally, thanks must also go to the anonymous referee, whose suggestions and careful reading of the manuscript led to many further improvements to the paper.
The author is funded by the internal graduate studentship of Trinity College, Cambridge.






\bibliographystyle{acm}  
\renewcommand{\bibname}{Bibliography}
\bibliography{main}

\begin{thebibliography}{10}

\bibitem{abraham1987note}
{\sc Abraham, U.}
\newblock A note on {D}ilworth's theorem in the infinite case.
\newblock {\em Order 4\/} (1987), 107--125.

\bibitem{aharoni1984konig}
{\sc Aharoni, R.}
\newblock {K}{\"o}nig's duality theorem for infinite bipartite graphs.
\newblock {\em Journal of the London Mathematical Society 2}, 1 (1984), 1--12.

\bibitem{aharoni1987menger}
{\sc Aharoni, R.}
\newblock Menger's theorem for countable graphs.
\newblock {\em Journal of Combinatorial Theory, Series B 43}, 3 (1987), 303--313.

\bibitem{aharoni1991infinite}
{\sc Aharoni, R.}
\newblock Infinite matching theory.
\newblock {\em Discrete Mathematics 95}, 1-3 (1991), 5--22.

\bibitem{aharoni2022strongly}
{\sc Aharoni, R.}
\newblock Strongly maximal matchings and strongly minimal covers.
\newblock {\em arXiv preprint arXiv:2206.02576\/} (2022).
\newblock 3 pages.

\bibitem{aharoni2009menger}
{\sc Aharoni, R., and Berger, E.}
\newblock Menger’s theorem for infinite graphs.
\newblock {\em Inventiones mathematicae 176}, 1 (2009), 1--62.

\bibitem{aharoni2011strongly}
{\sc Aharoni, R., and Berger, E.}
\newblock Strongly maximal antichains in posets.
\newblock {\em Discrete mathematics 311}, 15 (2011), 1518--1522.

\bibitem{aharoni1994menger}
{\sc Aharoni, R., and Diestel, R.}
\newblock Menger's theorem for a countable source set.
\newblock {\em Combinatorics, Probability and Computing 3}, 2 (1994), 145--156.

\bibitem{aharoni1992greene}
{\sc Aharoni, R., and Korman, V.}
\newblock Greene-{K}leitman's theorem for infinite posets.
\newblock {\em Order 9\/} (1992), 245--253.

\bibitem{aharoni1995strongly}
{\sc Aharoni, R., and Loebl, M.}
\newblock Strongly perfect infinite graphs.
\newblock {\em Israel Journal of Mathematics 90}, 1 (1995), 81--91.

\bibitem{diestel2003countable}
{\sc Diestel, R.}
\newblock The countable {E}rd{\H{o}}s--{M}enger conjecture with ends.
\newblock {\em Journal of Combinatorial Theory, Series B 87}, 1 (2003), 145--161.

\bibitem{dilworth1950decomposition}
{\sc Dilworth, R.}
\newblock A decomposition theorem for partially ordered sets.
\newblock {\em Annals of Mathematics 51}, 1 (1950), 161--166.

\bibitem{duffus2002intervals}
{\sc Duffus, D., and Goddard, T.}
\newblock Some progress on the {A}haroni--{K}orman conjecture.
\newblock {\em Discrete mathematics 250}, 1-3 (2002), 79--91.

\bibitem{goddard1996ordered}
{\sc Goddard, T.}
\newblock {\em Ordered sets: Colorings and complexity}.
\newblock {PhD Thesis}, Emory University, June 1996.

\bibitem{greene1976structure}
{\sc Greene, C., and Kleitman, D.~J.}
\newblock The structure of {S}perner k-families.
\newblock {\em Journal of Combinatorial Theory, Series A 20}, 1 (1976), 41--68.

\bibitem{hausdorff1908classification}
{\sc Hausdorff, F.}
\newblock Grundz{\"u}ge einer {T}heorie der geordneten {M}engen.
\newblock {\em Mathematische Annalen 65}, 4 (1908), 435--505.

\bibitem{hollom2024resolution}
{\sc Hollom, L.}
\newblock A resolution of the aharoni-korman conjecture.
\newblock {\em arXiv preprint arXiv:2411.16844\/} (2024).

\bibitem{mehta2025formal}
{\sc Mehta, B.}
\newblock {Disproof of the Aharoni-Korman conjecture}.
\newblock \url{https://github.com/b-mehta/AharoniKorman}, 2025.
\newblock A formal verification of the proof of Proposition 5.7. Accessed 2025-01-15.

\bibitem{perles1963dilworth}
{\sc Perles, M.~A.}
\newblock On {D}ilworth’s theorem in the infinite case.
\newblock {\em Israel Journal of Mathematics 1\/} (1963), 108--109.

\bibitem{van2022counterexample}
{\sc van~der Zypen, D.}
\newblock Counterexample to a conjecture of {A}haroni and {K}orman.
\newblock {\em arXiv preprint arXiv:2205.02296\/} (2022).
\newblock 2 pages.

\bibitem{zaguia2024progress}
{\sc Zaguia, I.}
\newblock {Some progress on the Aharoni–Korman conjecture}.
\newblock {\em Discrete Mathematics 347}, 10 (2024).

\end{thebibliography}


\appendix




\end{document}